\documentclass[preprint,12pt]{elsarticle}

\usepackage{amssymb}

\newtheorem{assumption}{Assumption}
\newtheorem{lemma}{Lemma}
\newtheorem{theorem}{Theorem}
\newtheorem{corollary}{Corollary}
\newtheorem{remark}{Remark}

\def\be{\begin{equation}}
\def\ee{\end{equation}}
\def\bea{\begin{eqnarray}}
\def\eea{\end{eqnarray}}
\newcommand{\re}[1]{(\ref{#1})}
\def\ne{ \nonumber \\ }
\def\nn{ \nonumber }

\def\QED{~\rule[-1pt]{5pt}{5pt}\par\medskip}

\begin{document}

\begin{frontmatter}



\title{Convergence of incremental adaptive systems}


\author{Mingxuan Sun}

\address{College of Information Engineering, Zhejiang University of Technology\\ Hangzhou 310023, CHINA}

\begin{abstract}
In this paper, incremental adaptive mechanisms
are presented and characterized,
to provide design hints for the development of continuous-time adaptive systems.
The comparison with the conventional integral adaptive systems indicates that the suggested design methodology
will be a supplement to the existing ones.
With the aid of a Barbalat-like lemma,
convergence results of the incremental adaptive systems are established.
It is shown that the proposed adaptive mechanisms are able to work well
in handling parametric uncertainties
in systems undertaken.
\end{abstract}


\end{frontmatter}


\section{Introduction}
\label{sec.intro}

Consider the following uncertain system
\bea
\dot x &=& \theta^{0 T} \varphi(t,x)  + b u
\label{sys}
\eea
where $x$ is the scalar state, $u$ is the scalar input of the system,
$\theta^{0}$ is the vector of unknown parameters,
$\varphi(\cdot,\cdot)$ is the vector of known continuous nonlinearities, and
$b$ is the unknown control gain,
but its sign is assumed to be known.
Here,
we does not specify whether $b$ is positive or negative,
implying that the adaptive systems to be developed in this paper are suitable to both situations.

It is seen that the parameter vector $\theta^{0}$ appear linearly in \re{sys}, and this indicates
that the main point of this paper is to handle the linear-in-the-parameters uncertainty.
The problem is, for a given trajectory $x_d(t), t \in [0, +\infty),$
to develop adaptive mechanisms
for estimating the unknown parameters,
and based on the estimates find the control   $u(t), t \in [0, +\infty),$
such that $x(t)$ follows $x_d(t)$ as close as possible, as $t \rightarrow +\infty.$

Let us denote by $e = x - x_d$ the tracking error.
The time derivative of the tracking error with respect to time can be expressed as
\bea
\dot e &=& b(  \theta^T \varphi(t,x)  + u)
\label{e.theta}
\eea
with $\theta = \theta^0/b$.

Let us begin with a discussion on the conventional integral adaptive systems.
We refer the reader to literature
\cite{goodwin87,middleton88,sastry89},
for design issues in continuous-time adaptive control,
\cite{narendra89}
for model reference adaptive methodologies,
\cite{wen95} for
robustness of adaptive systems,
\cite{ioannou96,tao03} for robust adaptive algorithms,
\cite{krstic95} for
adaptive backstepping designs,
and \cite{astolfi07}
more recent
immersion and invariance adaptive techniques.
As is well known, it is difficult to establish the asymptotic stability of time-varying systems
as it is not easy to find the Lyapunov function with a negative definite derivative.
Fortunately,
Barbalat's Lemma is found to be useful in performance analysis of adaptive systems,  which states that
if the integral of a uniformly continuous function having a (finite) limit, then the function
converges to zero asymptotically\cite{krstic95,ioannou96}.
Note that a simple alternative to Barbalat's Lemma can be found in
\cite{tao03}.

Consider the adaptive system consisting of system \re{sys}, the controller
 \bea u = - {\rm sgn} (b) \kappa e - \hat \theta^T \varphi \label{u.ac}
 \eea
and the adaptation law
 \bea
 \dot {\hat\theta} = - {\rm sgn} (b) \gamma \varphi e
  \label{al.ac}
 \eea
where $\kappa, \gamma > 0$ are design parameters, and ${\rm sgn} (\cdot)$ is the sign function.
Controller \re{u.ac} is designed based on the certainty equivalent principle.
To establish convergence of the adaptive system, we choose the positive definite function,
$L = V + \frac{|b|}{2\gamma } \tilde \theta^T\tilde \theta$, where $V = \frac{1}{2} e^2$,
and $\tilde \theta = \hat \theta - \theta$.
Applying \re{u.ac} and \re{al.ac},
we have
$ \dot L  = -  c_1 V $, $c_1 = 2\kappa |b| $.
This implies that $\dot L$ is
negative semidefinite, which renders $L$ to be bounded.
Due to the boundedness of $L$, it is easy to obtain
the boundedness of $V$ and its derivative,
 $\dot V$,
as well as
$\int_0^t V(s) ds < +\infty $,
as $t \rightarrow +\infty$.
Invoking Barbalat's Lemma shows that $\lim_{t \rightarrow + \infty} V(t) = 0$.
In turn, we conclude that
$
\lim_{t \rightarrow + \infty} e(t) = 0
$.

Now let us look into Eq. \re{al.ac}, by integrating its both sides as follows:
\bea
 \hat{\theta}(t) = \hat{\theta}(0) - {\rm sgn} (b) \gamma \int_{0}^t \varphi(s,x(s)) e(s) ds
 \label{eq.al.integral.1}
 \eea
which give the estimate for $\theta$
through the indicated integration.
The adaptation law \re{al.ac} is usually referred to as an integral adaptive law.
As for $t> \tau$,
\bea
 \hat{\theta}(t- \tau) = \hat{\theta}(0) - {\rm sgn} (b) \gamma  \int_{0}^{t- \tau}  \varphi(s,x(s)) e(s) ds
 \label{eq.al.integral.2}
 \eea
Subtracting \re{eq.al.integral.2} from \re{eq.al.integral.1}, we obtain
\bea
 \hat{\theta}(t) = \hat{\theta}(t-\tau) - {\rm sgn} (b) \gamma \int_{t-\tau}^{t} \varphi(s,x(s)) e(s) ds, ~~ t > \tau \nn
 \eea
Then appealing to the integral mean-value theorem, an incremental form of \re{al.ac} is obtained as follows:
\bea
\hat{\theta}(t) = \hat{\theta}(t-\tau) - {\rm sgn} (b) \tau \gamma \varphi(t_1, x(t_1))  e(t_1), ~~ t >\tau
 \label{eq.al.integral.4}
 \eea
where $t_1$ lies between $t-\tau$ and $t$,
takes different values for different instants of time.
To use $x(t)$ or $x(t-\tau)$ are two ways
to approximate $x(t_1)$
for implementation of the incremental adaptive mechanism.

Transient performance is always a major concern
in an adaptive system design.
Due to slow rate of parameter convergence,
it may exhibit poor transient behavior together with ideal asymptotic performance.
We are concerned about the second term of the right-hand side of \re{eq.al.integral.4}, where $\tau$ appears.
The parameter estimates will become hard to adapt,
as $\tau$ is set to be small.
One way is to reduce the sampling rate.
However,
the parameter estimates would not in time updated with the measured data, when setting $\tau$ too large.
As such,
the expression of \re{eq.al.integral.4} suggests that we choose $\gamma$ to be proportional to $\tau$ as follows:
\bea
 \gamma = \frac{1}{\tau} \gamma'
 \label{eq.Gamma}
 \eea
with $\gamma'$ being a constant specified by designer.
This discussion is an motivation for this paper to suggest the novel methodology of incremental adaptation.
Unlike the conventional ones,
the incremental adaptive mechanisms
do not involve such a $\tau$.
We shall clarify in the next section
how the incremental adaptive mechanisms are different from the conventional ones.

\section{Analysis of incremental adaptive systems}

Barbalat's lemma is a tool commonly employed for
concluding the convergence results of integral adaptive systems.
The following presents a Barbalat-like lemma,
a slightly modified form of
Lemma 1 in \cite{sun09} and
Lemma 1 in \cite{sun12},
which is specifically tailored for analysis of the incremental adaptive systems.

\begin{lemma} \label{lem.barbalatlike}
Suppose that $g(t)$, a time function positive on $[0, +\infty)$, satisfies
 \bea \int_{t-\tau}^t \dot g^2(s) ds \le M
 \label{eq.g2}
 \eea
for $t \in [\tau, +\infty)$, with $\tau >0$ being a constant, and
\bea
\lim_{t \rightarrow + \infty} \int_{t-\tau}^t g(s) ds =0
\label{eq.intg}
\eea
Then $\lim_{t \rightarrow + \infty} g(t) =0$.
\end{lemma}
{\bf Proof.}
See Appendix for the proof.
\QED

\begin{corollary}  \label{cor.barbalat}
Lemma \ref{lem.barbalatlike} holds, if the condition \re{eq.g2} is replaced with
that $\dot g(t)$ is bounded.
\end{corollary}
{\bf Proof.} This corollary follows immediately from the observation that Eq. \re{eq.g2} holds, whenever $\dot g(t)$ is bounded.
\QED

We are now in a position to present
the convergence result of an adaptive system, where
the incremental adaptation mechanism is adopted.

\begin{theorem} \label{thm.iac}
Consider the incremental adaptive system described by the system \re{sys}, the controller
 \bea u = - {\rm sgn} (b) \kappa e - \hat \theta^T \varphi \label{u.iac}
 \eea
and the adaptation law
  \bea
\hat\theta(t) =
\left \{ \begin{array}{ll}
\hat\theta(t-\tau) - {\rm sgn} (b) \gamma \varphi(t,x(t)) e(t) & ~~ for ~~ t \ge 0  \\
 \hat \theta_0 & ~~ for ~~ t \in [-\tau, 0)
  \end{array}
 \right.
  \label{al.iac}
 \eea
where $\tau, \kappa, \gamma > 0$ are design parameters, and $ \hat \theta_0$ is the initial setting for $\hat\theta$.
Then the tracking error $e(t)$ will be guaranteed to converge to zero, as time increases, i.e.,
\bea
\lim_{t \rightarrow + \infty} e(t) = 0 \nn
\eea
while
$e, x$, as well as $\int_{t-\tau}^{t} \|\hat \theta(s)\|^2 ds$ and
  $ \int_{t-\tau}^{t} u^2(s) ds$,
for $t \in [\tau, + \infty),$
are bounded.
\end{theorem}

{\bf Proof.} The convergence result of the closed-loop system
composed of \re{sys}, \re{u.iac} and \re{al.iac} can be established by choosing
the following Lyapunov-Krasovskii functional candidate
 \bea
 L(t) = V(t) + \frac{|b|}{2\gamma } \int_{t-\tau}^t \tilde \theta^T(s) \tilde \theta(s) ds
 \label{L.iac}
 \eea
with $V = \frac{1}{2} e^2$. Its derivative with respect to time is
\bea
 \dot L(t)
 &=& \dot V(t) + \frac{|b|}{2\gamma } [ \tilde \theta^T(t) \tilde \theta(t)
- \tilde \theta^T(t-\tau)\tilde \theta(t-\tau)
 ]  \label{eq.dL}
 \eea
To proceed, the filtered error dynamics is expressed as, when applying \re{u.iac},
\bea
\dot e
&=& - \kappa |b| e
+ b \tilde \theta^T \varphi \label{eq.ef}
\eea
The derivative of $V$ along the error trajectory \re{eq.ef} can be given as
\bea
 \dot V
 &=& - \kappa |b| e^2
+ b \tilde \theta^T \varphi e
\label{eq.dV.1}
 \eea
The second term of the right-hand side of \re{eq.dL} satisfies
 \bea
 &&\tilde \theta^T(t-\tau) \tilde \theta(t-\tau) \ne
 &=& [\hat \theta(t-\tau)- \hat\theta(t)
 + \hat\theta(t) - \theta]^T
 [\hat \theta(t-\tau)- \hat\theta(t)
 + \hat\theta(t) - \theta] \ne
 &=& [\hat\theta(t)-\hat\theta(t-\tau)]^T[\hat\theta(t)-\hat\theta(t-\tau)] \ne
 &&+ 2[\hat \theta(t) - \hat\theta(t-\tau)]^T \tilde{\theta}(t)
  + \tilde \theta^T(t) \tilde \theta(t) \label{eq.theta.1}
 \eea
Substituting \re{eq.dV.1} and \re{eq.theta.1} into \re{eq.dL}, we obtain
\bea
 \dot L(t)
  &=& -\kappa |b| e^2(t)
+ b \tilde \theta^T(t) \varphi(t,x(t)) e(t) \ne &&
   + \frac{|b|}{\gamma} \tilde{\theta}^T(t) [\hat\theta(t)-\hat\theta(t-\tau)] \ne &&
   - \frac{|b|}{2\gamma }
[\hat\theta(t)-\hat\theta(t-\tau)]^T
[\hat\theta(t)-\hat\theta(t-\tau)]
 \nn
 \eea
Then applying \re{al.iac} yields
\bea
 \dot L(t)
  &=& -\kappa |b| e^2(t)
   - \frac{|b|}{2\gamma }[\hat\theta(t)-\hat\theta(t-\tau)]^T
[\hat\theta(t)-\hat\theta(t-\tau)]\ne
  &\le&  - \kappa |b| e^2(t) \ne
    &=& - c_1 V(t)
  \label{eq.dL.1}
 \eea
where $c_1 = 2 \kappa |b|$.

Eq.\re{eq.dL.1} makes $\dot L$  negative semidefinite.
The boundedness of $L$ is ensured due to the boundedness of $L(0)$.
Hence, $V$ is bounded, implying the boundedness of $e$, and in turn that of $x$.
Furthermore,
$\int_{t-\tau}^t \|\hat \theta(s)\|^2 ds$,
$t \in [\tau, + \infty),$
 is bounded, by the definition of $L$.
It follows from \re{u.iac} that
 \bea u^2 &\le& 2 \kappa^2 e^2 + 2 (\hat \theta^T \varphi )^2 \ne
  &\le& c_2 + c_3 \|\hat \theta\|^2
 \label{u2.iac}
 \eea
where
$c_2 =   2\kappa^2 \sup_{t \in [0,+ \infty)} e^2$, and
$c_3 = 2 \sup_{t \in [0,+ \infty)} \|\varphi\|^2$. Hence,
the boundedness of $\int_{t-\tau}^t u^2(s)ds$, $t \in [\tau, + \infty),$  follows by noting that
\bea \int_{t-\tau}^t u^2(s)ds &\le& c_2 \tau + c_3
\int_{t-\tau}^{t} \|\hat \theta(s)\|^2 ds < + \infty
 \label{u2.iac.1}
 \eea

The difference between instants of $t$ and $t-\tau $,  $L(t) - L(t-\tau)$, can be calculated by
\bea
 L(t) - L(t-\tau)
 &=& \int_{t-\tau}^t \dot L(s) ds
 \nn
 \eea
Again using \re{eq.dL.1},
 \bea L(t) - L(t-\tau) \le - c_1 \int_{t-\tau}^t V(s) ds
 \nn
 \eea
For each fixed instant $t = t_i = i \tau + t_0$,  $i = 1,2,\dots$, $t_0 \in [0, \tau)$.
 \bea
 L(t_i) - L(t_{i-1}) &\le& - c_1  \int_{t_{j-1}}^{t_j} V(s) ds
 \nn
 \eea
leading to
 \bea
 L(t_i) - L(t_0) &\le& -c_1\sum_{j=1}^i \int_{t_{j-1}}^{t_j} V(s) ds
 \nn
 \eea
Consequently, by the finiteness of $L(t_0)$, the series $\sum_{j=1}^i \int_{t_{j-1}}^{t_j} V(s) ds$ converges. Therefore,
\bea \lim_{i \rightarrow \infty } \int_{t_{i-1}}^{t_i} V(s) ds
&=& 0
 \nn
 \eea
implying that
\bea
 \lim_{t \rightarrow \infty } \int_{t-\tau}^{t}
V(s) ds &=& 0 \label{eq.intV}
 \eea

Now, we consider the finiteness of $\int_{t-\tau}^t\dot V^2(s)ds$, $t \in [\tau, + \infty)$.
It follows from \re{eq.dV.1} that
\bea
 \dot V ^2
 &=& [ - \kappa |b| e^2
+ b \tilde \theta^T \varphi e  ]^2 \ne
 &\le& 2 \kappa^2 b^2 e^4 + 2 b^2 (\tilde \theta^T \varphi)^2 e^2 \ne
 &\le& c_4 + c_5 \|\tilde \theta\|^2
   \nn
 \eea
where
$c_4 = 2\kappa^2 b^2 \sup_{t \in [0, +\infty)} e^4$ and $
c_5 = 2 b^2 \sup_{t \in [0, +\infty)}\|\varphi\|^2 \sup_{t \in [0, +\infty)}e^2$.
Integrating both sides and by the boundedness of $\tilde \theta$ yield
\bea
 \int_{t-\tau}^t \dot{V}^2(s) ds
 &\le& c_4 \tau + c_5 \int_{t-\tau}^t\|\tilde \theta(s)\|^2 ds < + \infty
 \label{eq.dV2}
 \eea
for $t \in [\tau, + \infty)$.
In view of \re{eq.intV} and \re{eq.dV2},
and by Lemma \ref{lem.barbalatlike},
$
\lim_{t \rightarrow + \infty} V(t) = 0.
$
In turn, we can conclude the convergence of $e(t)$, as $t \rightarrow + \infty$.
\QED

\begin{remark}
In comparison with \re{eq.al.integral.4},
no $\tau$ appears in the second term of the right-hand side of
the adaptation law \re{al.iac},
which indicates the main difference between
the integral adaptation law and the incremental
adaptation law.
\end{remark}

\begin{remark}
The adaptation law given in Theorem \ref{thm.iac} guarantees the boundedness of $\hat\theta$
in the sense as presented.
In order to ensure the boundedness of $\hat\theta$ itself,
the saturated learning is helpful \cite{sun06}.
In particular,
for fully saturated learning, the entire
right-hand side of the learning law is saturated,
and the estimate is ensured to be within a
pre-specified region.
We apply the fully-saturated adaptation law as follows:
  \bea
\left \{\begin{array}{l}
\hat\theta(t) = {\rm sat} (\hat\theta^*(t)) \\
\hat\theta^*(t) = {\rm sat} (\hat\theta^*(t-\tau)) - {\rm sgn} (b) \gamma \varphi(t,x(t)) e(t)
  \end{array}
 \right.
  \label{al.sat.iac}
 \eea
$for ~~ t \ge 0$.
By the boundedness of $e$, $x$ and $\hat\theta$,
it is easy to obtain the
boundedness of
$u$ from \re{u.ac}, and that of $\dot V$ from \re{eq.dV.1}. By invoking Corollary \ref{cor.barbalat}, convergence of such an incremental adaptive system can be established.
\end{remark}

The following theorem clarifies the flexibility of choice of incremental adaptive mechanisms.

\begin{theorem} \label{thm.iac.openloop}
When the adaptive control law given by
 \bea u &=& - {\rm sgn} (b) \kappa e - \hat \theta^T \varphi + u_1 \label{u.iac.open} \\
 u_1 &=&  - \frac{1}{2}  {\rm sgn} (b) \gamma\varphi^2 e \label{u1.iac}
 \eea
with the adaptation law
  \bea
\hat\theta(t + \tau) =
\left \{ \begin{array}{ll}
\hat\theta(t) - {\rm sgn} (b) \gamma \varphi(t,x(t)) e(t) & ~~ for ~~ t \ge 0  \\
 \hat \theta_0 & ~~ for ~~ t \in [0, \tau)
  \end{array}
 \right.
  \label{al.iac.open}
 \eea
is applied to system \re{sys},
then the same results as in Theorem \ref{thm.iac} are true.
\end{theorem}

{\bf Proof.}
In order to cope with
the use of \re{al.iac.open},
we choose the following candidate Lyapunov-Krasovskii functional,
 \[
 L(t) = V(t) + \frac{|b|}{2\gamma } \int_t^{t+\tau} \tilde \theta^T(s) \tilde \theta(s) ds
 \]
with the same $V(t)$ as that in \re{L.iac}.
Employing \re{u.iac.open}, the error dynamics can be expressed as
\bea
\dot e
&=& - \kappa |b| e
+ b \tilde \theta^T \varphi + b u_1 \nn
\eea
The derivative of $V$ along trajectories of the error dynamics is given by
\bea
 \dot V
 &=& e \dot e \ne
 &=& - \kappa |b| e^2
+ b \tilde \theta^T \varphi e + b u_1 e
 \nn
 \eea
Hence, the derivative of $L$ can be calculated as
\bea
 \dot L(t)
 &=& \dot V(t) + \frac{|b|}{2\gamma } [ \tilde \theta^T(t+\tau)\tilde \theta(t+\tau) - \tilde \theta^T(t) \tilde \theta(t)
 ]  \nn
 \eea
The second term of the right-hand side of
the above equation satisfies
 \bea
 &&\tilde \theta^T(t+\tau) \tilde \theta(t+\tau) \ne
 &=& [\hat \theta(t+\tau)- \hat\theta(t)
 + \hat\theta(t) - \theta]^T
 [\hat \theta(t+\tau)- \hat\theta(t)
 + \hat\theta(t) - \theta] \ne
 &=& [\hat\theta(t+\tau)-\hat\theta(t)]^T
 [\hat\theta(t+\tau)-\hat\theta(t)] \ne
 &&+ 2[\hat \theta(t+\tau) - \hat\theta(t)]^T \tilde{\theta}(t)
  + \tilde \theta^T(t)\tilde \theta(t) \nn
 \eea
It follows that
\bea
 \dot L(t)
  &=& -\kappa |b| e^2(t)
+ b \tilde \theta^T(t) \varphi(t,x(t)) e(t)  + b u_1(t) e(t)\ne &&
   + \frac{|b|}{\gamma} \tilde{\theta}^T(t) [\hat\theta(t+\tau) - \hat\theta(t)] \ne &&
   + \frac{|b|}{2\gamma }
[\hat\theta(t+\tau) - \hat\theta(t)]^T
[\hat\theta(t+\tau) - \hat\theta(t)]
 \nn
 \eea
Applying \re{u1.iac} and \re{al.iac.open} yields
\bea
\dot L &=& -\kappa |b| e^2 \nn
 \eea
The proof can be carried out
by evaluating the term
$ L(t)-L(t-\tau)$,
with similar lines to those of the proof for Theorem \ref{thm.iac}.
\QED

\begin{remark}
It is seen in \re{u.iac.open} that
an additional component, $u_1$, is added into \re{u.iac},
for canceling the term appeared when applying
\re{al.iac.open}.
\end{remark}

\section{Robust treatments}
In this section,
we shall provide an approach for analysis of
the adaptive system to be developed, in the presence of bounded uncertainty, by
considering the class of single-input single-output continuous-time systems
\bea
y^{(n)} + \sum_{i=1}^{n_a} a_i Y_i(t, y, \dot y,\cdots, y^{(n-1)} )  &=& b u + w
\label{sys.y}
\eea
where $u$ and $y$ are the scalar input and output of the system, respectively, and $w$ represents the lumped non-parametric uncertain term;
$a_i, i=1,2, \cdots, n_a,$ are unknown coefficients, and
$Y_i, i=1,2, \cdots, n_a, $ represent known nonlinearities, being bounded as
$y, \dot y, \cdots, y^{(n-1)}$ are bounded; and
$b$ is the unknown control gain.

By introducing the state vector $x = [x_1, \cdots, x_n]^T$,
and the state space representation for system \re{sys.y} can be given as follows:
 \bea
\left \{ \begin{array}{l}
 \dot x_i  = x_{i+1}, i = 1, 2, \cdots, n-1\\
 \dot x_n = - \sum_{i=1}^{n_a} a_i Y_i(t, x) + b u + w \\
  y = x_1
 \end{array}
 \right.
 \label{sys.x}
 \eea

Let us denote by $e = x - x_d = [e_1, e_2, \cdots, e_n]^T$ the tracking error, where
$x_d = [y_d, \dot y_d, \cdots, y_d^{(n-1)} ]^T,$ and for $\lambda > 0$,
$e_f = \left ( \frac{d}{dt} + \lambda \right)^{n-1} e_1$ the filtered error, where
$y_d(t), t \in [0, +\infty),$ is the  desired trajectory.
The time derivative of $e_f$ with respect to time is of the form
\bea
\dot e_f &=&  - \sum_{i=1}^{n_a} a_i Y_i(t, x) + b u + w + \nu
\label{ef}
\eea
with
$\nu = [0 ~~\Lambda^T] e - y_d^{(n)}$ and $\Lambda =[\lambda^{n-1}, (n-1)\lambda^{n-2} , \cdots, (n-1)\lambda ]^T $.

\begin{assumption} \label{assump.b}
The sign of the control gain $b$ is known.
\end{assumption}

As discussed before,
we again does not specify whether $b$ is positive or negative.
Define $\theta = [a_1/b, \cdots, a_{n_a}/b, 1/b]^T$ and $\varphi(t,x) = [-Y_1(t, x), \cdots,$ $ -Y_{n_a}(t,x), \nu]^T$. Eq. \re{ef} can be rewritten as
\bea
\dot e_f &=& b(  \theta^T \varphi(t,x)  + u + w_b)
\label{ef.theta}
\eea
where $w_b = w/b$.

\begin{assumption} \label{assump.wb}
The uncertain term $w_b$ is assumed to be bounded, satisfying
\bea
|w_b| &\le & \bar w_b
\label{wb}
\eea
where $\bar w_b = \frac{\bar w}{|b|},$ and $|w| \le \bar w$.
\end{assumption}

Now we present the robust treatments in forming an incremental adaptation mechanism in the presence of $w_b$.
Let us introduce functions $\iota_\epsilon(\cdot)$ and $\varsigma_\epsilon(\cdot)$ as follows:
\bea
 \iota_\epsilon (\cdot) = \left\{ \begin{array}{ll}
    1
  & \mbox{if $|\cdot| > \epsilon$}  \\
     0 & \mbox{if $|\cdot| \le \epsilon$}
   \end{array}
   \right.
\label{def.iota}
 \eea
and
\bea
 \varsigma_\epsilon (\cdot) = \left\{ \begin{array}{ll}
       {\rm sgn}(\cdot)
     & \mbox{if $|\cdot| > \epsilon$}  \\
   0 & \mbox{if $|\cdot| \le \epsilon$}
    \end{array}
   \right.
\label{def.varsigma}
 \eea
and define the error variable
$e_{\epsilon}(t)
= (|e_f(t)| - \epsilon)\iota_{\epsilon}(t)$.

\begin{theorem} \label{thm.riac}
Consider the incremental adaptive system described by the system \re{sys.y}, the controller
 \bea u =
 - {\rm sgn} (b) \kappa e_\epsilon \varsigma_\epsilon
  - {\rm sgn} (b) \bar w_b \iota_\epsilon
 - \hat \theta^T \varphi \iota_\epsilon
\label{u.riac}
 \eea
and the adaptation law
  \bea
\hat\theta(t) =
\left \{ \begin{array}{ll}
\hat\theta(t-\tau)
- {\rm sgn}(b) \gamma
\varphi(t,x(t)) e_{\epsilon}(t)
\varsigma_{\epsilon}(t)
 & ~~ for ~~ t \ge 0  \\
 \hat \theta_0 & ~~ for ~~ t \in [-\tau, 0)
  \end{array}
 \right.
  \label{al.riac}
 \eea
where $\tau, \kappa, \gamma > 0$ are parameters to be specified by designer, and $ \hat \theta_0$ is the initial setting for $\hat\theta$.
Then the error variable $e_{\epsilon}(t)$ can be made to converge to zero, as time increases, i.e.,
\bea
\lim_{t \rightarrow + \infty} e_{\epsilon}(t) = 0 \nn
\eea
while
$e_{\epsilon}, e_f, e, x,$ as well as $ \int_{t-\tau}^{t} \|\hat \theta(s)\|^2 ds$ and
$ \int_{t-\tau}^{t} u^2(s) ds$, $t \in [\tau, + \infty)$,
are all bounded.
\end{theorem}

{\bf Proof.}
The proof follows similar lines to those of the proof of Theorem \ref{thm.iac},
with the positive definite function
\bea V = \frac{1}{2} e_{\epsilon}^2
 \label{v.riac}
\eea
By \re{def.iota} and \re{def.varsigma}, the derivative of $V$ with respect to time is calculated as
 \bea
 \dot{V}
   &=& e_{\epsilon} \varsigma_\epsilon \dot e_f \ne
  &=& e_{\epsilon} \varsigma_\epsilon b [ \theta^T \varphi + u  + w_b]\ne
 &\le&
 - |b| \kappa e_\epsilon^2
 -  b \tilde \theta^T \varphi e_{\epsilon} \varsigma_\epsilon
   \label{dv.riac}
 \eea
We take the same positive definite function \re{L.iac}
as a Lyapunov-Krasovskii functional candidate,
with the defined $V(t)$ in \re{v.riac}. The derivative of $L(t)$ can be calculated as
 \bea
 \dot L(t) &=&
 \dot V(t) +
  \frac{|b|}{2\gamma} (\tilde \theta^T(t) \tilde \theta(t)
- \tilde \theta^T(t-\tau) \tilde \theta(t-\tau) )
 \label{dL.riac.1}
 \eea
Using \re{eq.theta.1} and \re{dv.riac},
$\dot L(t)$ given by Eq. \re{dL.riac.1}
satisfies
\bea
 \dot L(t) &\le&
 - |b|\kappa e_{\epsilon}^2(t) - b \tilde \theta^T(t) \varphi(x(t),t) e_{\epsilon}(t) \varsigma_\epsilon(t) -
  \frac{|b|}{\gamma} \tilde \theta^T(t) [\hat{\theta}(t) -
\hat{\theta}(t-\tau)]  \ne
  &&-
  \frac{|b|}{2\gamma}
[\hat{\theta}(t) -
\hat{\theta}(t-\tau)]
^T [\hat{\theta}(t) -
\hat{\theta}(t-\tau)]
  \nn
 \eea
Applying the adaptation law \re{al.riac}, we obtain
\bea
 \dot L(t) &\le&
 - |b|\kappa e_{\epsilon}^2(t) - \frac{1}{2\gamma}
[\hat{\theta}(t) -
\hat{\theta}(t-\tau)]^T [\hat{\theta}(t) -
\hat{\theta}(t-\tau)]
  \ne
&\le&
  - |b|\kappa e_{\epsilon}^2(t) \ne
&=& - c_1 V(t)
  \label{dL.riac}
 \eea
where $c_1 = 2|b|\kappa$.

By \re{dL.riac},
the boundedness of $L$ is ensured, as $\dot L$ is negative semidefinite and $L(0)$ is bounded.
Thus, by the definition of $L$, $V$ is bounded, implying the boundedness of $e_{\epsilon}$, $e_f$, $e$, and in turn that of $x$.
Moreover, by the definition of $L$,
$\int_{t-\tau}^t \|\hat \theta(s)\|^2 ds$ is bounded for $t \in [\tau, +\infty)$.
With the similar derivations to those of the proof for Theorem \ref{thm.iac}, we can conclude
the boundedness of $\int_{t-\tau}^t u^2(s)ds$ for $t \in [\tau, +\infty)$.

To proceed, we recall the expression of the difference of $ L(t) - L(t-\tau)$, given by
\bea
 L(t) - L(t-\tau) &=&
 V(t) - V(t-\tau) \ne &&+
  \frac{|b|}{2\gamma} \int_{t-\tau}^t (\tilde \theta^T(s) \tilde \theta(s)
- \tilde \theta^T(s-\tau) \tilde \theta(s-\tau) )ds \nn
 \eea
It follows from \re{dv.riac} that
\bea
 V(t) - V(t-\tau)
 &\le& - |b|\kappa \int_{t-\tau}^t  e_{\epsilon}^2(s) ds - b \int_{t-\tau}^t
\tilde \theta^T(s)\varphi(x(s),s) e_{\epsilon}(s) \varsigma_\epsilon(s) ds \nn
 \eea
and from \re{dL.riac},
\bea
 L(t)-L(t-\tau) &\le&
  - |b|\kappa \int_{t-\tau}^t  e_{\epsilon}^2(s) ds \ne
&& -
  \frac{1}{2\gamma}\int_{t-\tau}^t
[\hat{\theta}(s) -
\hat{\theta}(s-\tau)]
^T [\hat{\theta}(s) -
\hat{\theta}(s-\tau)]
  ds \ne
&\le&
  - |b|\kappa \int_{t-\tau}^t  e_{\epsilon}^2(s) ds \ne
&=&
  - c_1 \int_{t-\tau}^t  V(s) ds
  \nn
 \eea
where $ c_1 = 2|b|\kappa $. It follows that
 \bea
 L(t_i) - L(t_0) &\le& -c_1\sum_{j=1}^i \int_{t_{j-1}}^{t_j} V(s) ds
 \nn
 \eea
for each fixed instant $t = t_i = i \tau + t_0$,  $i = 1,2,\dots$, $t_0 \in [0, \tau)$.
Consequently, by the finiteness of $L(t_0)$,
\bea \lim_{i \rightarrow \infty } \int_{t_{i-1}}^{t_i} V(s) ds
&=& 0
 \nn
 \eea
implying that
\bea
 \lim_{t \rightarrow \infty } \int_{t-\tau}^{t}
V(s) ds &=& 0 \label{eq.intV.riac}
 \eea

Using \re{dv.riac} again,
 \bea
 \dot{V}^2
 &\le & 2[ b^2 \kappa^2 e_{\epsilon}^2 + [b \tilde \theta^T \varphi e_{\epsilon} \varsigma_\epsilon ]^2 ]\ne
  &\le & 2b^2 [ \kappa^2 + \|\tilde \theta\|^2 \| \varphi\|^2 ] e_{\epsilon}^2
   \nn
 \eea
leading to
\bea
\int_{t-\tau}^t \dot V ^2(s) ds < + \infty \label{intV2.riac}
\eea
for $t \in [\tau, + \infty)$,
which holds due to the boundedness of $\int_{t-\tau}^t \|\hat \theta(s)\|^2 ds$
for $t \in [\tau, + \infty)$.
In view of \re{eq.intV.riac} and \re{intV2.riac}, by using Lemma \ref{lem.barbalatlike}, we conclude that
$\lim_{t \rightarrow \infty} V(t) =
0,$
and in turn
$\lim_{t \rightarrow \infty} e_\epsilon(t) =
0$. This completes the proof.
\QED

\begin{remark}
Theorem \ref{thm.riac} indicates the convergence
of the error variable $e_\epsilon(t)$
of the incremental adaptive system, as time increases. In addition,
whenever $e_\epsilon(t)$ converging to zero,
$e(t)$ converges to the interval $(-\epsilon, \epsilon)$, as $t \rightarrow \infty.$
\end{remark}

\section{Concluding remarks}
We suggest
incremental adaptive mechanisms,
in this paper, applicable
to develop continuous-time adaptive systems, and illustrate
design hints for the development.
A comparison between
the integral and incremental adaptive systems is made to clarify why our approach makes sense.
It is interesting to note that
the update term of an incremental adaptation law looks the same as
the right-hand side term of the integral adaptation law, and
the integral adaptation law can be considered as a kind of incremental adaptive one,
where $\tau$, the duration of adaptation,
appears in the update term.
With the aid of the Barbalat-like lemma,
a unified approach for the analysis of incremental adaptive systems has been presented, by which the convergence has been established in the absence or presence of the disturbance term.
In this paper, we present our preliminary results on the
incremental adaptation.
For future work, we would like
to extend it to wide range of situations
where the conventional integral adaptive mechanisms are applicable.

\section*{Appendix}

For purpose of analysis, let us denote $t = k \tau + \sigma$, $\sigma \in [0, \tau)$, and $k = 0,1, \cdots$.
We prove by contradiction.
Suppose that we can find  $\sigma_0 \in [0,\tau)$, such that $g(k \tau + \sigma_0)$ does not converge to zero as $k \rightarrow \infty$. Then we know that there exist a subsequence $k_i$ and an $\varepsilon > 0$ such that
\bea
g(k_i \tau + \sigma_0) \ge \sqrt{\varepsilon}+
 \frac{1}{2} \sqrt{\sqrt{\varepsilon} M}  \label{a1}
\eea
and
\bea
\sigma_0 - \sqrt{\varepsilon}/2 \ge 0,~~
 \sigma_0 + \sqrt{\varepsilon}/2 < \tau \label{a2}
\eea

Select $\sigma_1 \in [\sigma_0 - \sqrt{\varepsilon}/2,
\sigma_0 + \sqrt{\varepsilon}/2] $.
By \re{eq.g2} and Schwarz's inequality, we obtain
\bea && |g( k_i \tau + \sigma_1) - g( k_i \tau + \sigma_0)| \ne &=&
\left | \int_{ k_i \tau + \sigma_0}^{ k_i \tau + \sigma_1} \dot g(s) ds \right | \ne
&\le&
\frac{1}{2} \int_{  k_i  \tau + \sigma_0 - \sqrt{\varepsilon}/2}^{  k_i  \tau + \sigma_0 + \sqrt{\varepsilon}/2} | \dot g (s) | ds  \ne
&\le& \frac{1}{2} \left (  \int_{  k_i  \tau + \sigma_0 - \sqrt{\varepsilon}/2}^{  k_i  \tau + \sigma_0 + \sqrt{\varepsilon}/2} 1^2 ds
 \int_{  k_i  \tau + \sigma_0 - \sqrt{\varepsilon}/2}^{  k_i  \tau + \sigma_0 + \sqrt{\varepsilon}/2}
|\dot g(s)|^2 ds \right)^{1/2} \ne
&\le& \frac{1}{2} \sqrt{ \sqrt{\varepsilon}M}  \label{a3}
\eea
Combining \re{a1} and \re{a3} yields
\bea |g(k_i \tau + \sigma_1) | & \ge & |g(k_i  \tau + \sigma_0) | - |g( k_i  \tau + \sigma_1) - g( k_i  \tau + \sigma_0)| \ne
& \ge & \sqrt{ \varepsilon}
 \label{a4}
\eea
It follows from \re{a4} that
\bea \int_{k_i \tau }^{(k_i+1) \tau } g(s)  ds  \ge \int_{k_i \tau + \sigma_0 - \sqrt{\varepsilon}/2
}^{k_i \tau + \sigma_0 + \sqrt{\varepsilon}/2 } g(s) ds  \ge
\sqrt{\varepsilon} \sqrt{\varepsilon}  = \varepsilon
\label{a5}
\eea
which contradicts to \re{eq.intg}. Therefore,
for each fixed $\sigma \in [0, \tau)$, $g(k \tau + \sigma)$ converges to zero as $k \rightarrow \infty$.
This completes the proof.




\end{document}